# A Set of Questions in Combinatorial and Metric Geometry

*R Nandakumar,*
*Amrita School of Arts and Sciences,*
*Edappalli North, Kochi – 682024, India.*
*([nandacumar@gmail.com](nandacumar@gmail.com))*

**Abstract:** We briefly introduce several problems: (1) a generalization of the 'convex fair partition' conjecture, (2) on non-trivial invariants among polyhedrons that can be formed from the same set of face polygons, (3) two questions on assembling rectangular tiles to form larger rectangles and (4) on convex regions which maximize and minimize the diameter for specified area and perimeter. For each question, we discuss partial solutions and indicate aspects that to our knowledge, await exploration.

## 1. On the 'Fair Partition' problem for convex regions

We attempt to generalize the fair partition conjecture ([1]) for convex regions and suggest that even if it is true, it may only be 'just true'.

**Definitions:** Given a planar convex shape, the *diameter* is defined to be the largest distance that can be formed between two parallel lines tangent to its boundary, and the *minimum width* is defined to be the smallest such distance.

**The original fair partition conjecture:** "Given a positive integer n, any convex 2D region allows partition into n convex pieces such that all pieces have the same area and perimeter." This conjecture has been established for prime power values of n and remains open for other values; a partial answer to the problem also has been recently named the 'spicy chicken theorem' ([13]).

**Generalization:** With n fixed at 2, let us generalize the question: Given 2 positive numbers a and b and if a convex region C is to be divided into 2 convex pieces with their areas in the ratio a : b and perimeters in the ratio $\sqrt{a} : \sqrt{b}$ (obviously, if a =b=1, we have the original convex fair partition problem). Such a partition, which considers the dimensions of the quantities, could be called a 'scaled fair partition'. This question could be generalized to higher dimensions: for example, in 3D if the convex pieces ought to have their diameters in the ratio a: b, the ratio between their surface areas is the ratio between their squares and the volume ratio is the ratio of their cubes. Generalization to higher n is obvious.

*A simple example:* Consider a rectangle of 1 X 4 units. Say, it is to be cut into 2 convex pieces with ratio between areas 1:3 and perimeters in ratio 1: $\sqrt{3}$ = 0.56
1. If the dividing line is parallel to the shorter side of the rectangle and the area ratio of the 2 resultant pieces is 1:3, and perimeter ratio is 4:8 = 0.5.
2. If the dividing line is parallel to the longer side of the rectangle and area ratio is 1:3, perimeter ratio is 8.5: 9.5 =0.89.

Between these two extreme orientations, the dividing line has a continuous range of orientations and as it is varied, the perimeter ratio will range continuously from 0.5 to 0.89. And this interval has the required perimeter ratio of 0.56.

However, we see that even if C is a circle, for piece areas in the ratio 1:3 say, such a scaled fair partition is not possible, the main reason being the requirement that the pieces be convex.

**Conjectures:** for general n, for any other ratio than a1: a2: .... : an = 1:1: ....:1 (the 'basic' fair partition case where all areas are equal and all perimeters are equal), there is always some convex shape which does not allow a scaled fair partition. Indeed, we suspect the circular disc might allow only the basic fair partition with all a's equal to 1 (Question: if so, is the disc the only such shape?) and hence the basic fair partition is somehow special for some deeper reason. In other words, even if the basic fair partition conjecture is correct (for non- prime power values of n, it is still open [13]), it may only be just correct.

**Non-convex pieces:** In version 2 of [2], we had proved that for any n, if the pieces are allowed to be non-convex, any convex shape can be divided into n pieces of same area and perimeter. Now, we sketch the proof for the following:

**Claim:** For n=2, with non-convex pieces allowed, any convex 2D region C allows scaled fair partition into 2 pieces, say P1 and P2, for any area ratio.

**Proof:** Indeed, if the required area ratio is a1 / a2 and the perimeter ratio is $\sqrt{(a1/a2)}$ and without loss of generality, a1 < a2, consider the following extreme configurations of P1 and P2: (1) the smaller piece P1 with area a1 is a convex and compact region with as small a perimeter as possible and shares a short portion of the boundary of the full region C; and the other piece P2 (obviously non-convex) is the rest of C. Here, the perimeter ratio P1: P2 is less than what is required for the partition into {P1, P2} to be scaled fair. (2) Both pieces P1 and P2 share the boundary of C equally and have the required areas (P1 would be a non-convex 'crescent-like' shape). Now, the perimeters of the two pieces P1 and P2 are equal and the perimeter ratio is 1, which is obviously greater than what is required. Now, since there is a continuous deformation of the separating curve between P1 and P2 that can take configuration 1 to configuration 2, there is an intermediate state where the perimeter ratio has the desired value.

For larger values of n and arbitrary area ratios, even with non-convex pieces, we don't know if a scaled fair partition always exists for any input convex region.

**Further questions:** In the basic convex fair partition, what are the bounds on the variation of diameter among the pieces with equal area and perimeter? [13] has proved the broad result: for n any prime power, a partition of any d-dimensional convex body into n convex pieces always exists that makes an absolutely continuous finite measure ( an example is volume) and any d-1 continuous functionals of the body (such as diameter, perimeter,…) equal among all the pieces. So one could ask (for example): in 2D (d=2), if the area and *more than 1* such continuous functionals need to be equal among the pieces, how would the possible convex partitions be restricted – for example the partition, if at all a partition exists, it could necessarily be a partition by fans? The idea is to see how the search space for partitions reduces. And what happens if the volumes need not be equal but more than d -1 continuous functionals need to be the same among n convex pieces?

# 2. On Reconstructing Convex Polyhedrons

**Basic question:** Consider a convex polyhedron. Break it at all edges until only face polygonal regions remain. Is it possible to do the reverse uniquely? i.e., given a set of disconnected polygons (of course, with some edges of these different 'faces' having equal lengths so they can be hinged along them), can one always get a unique convex polyhedron?

*Note 1:* the problem does not allow faces to be merged into larger faces. Else, given (say) a set of identical squares, one can merge them into various sets of 6 rectangles each and form many different cuboids.

*Note 2:* In 2D, the problem is trivial. Indeed, given a set of line segments, one could have convex polygons of hugely varying areas - with 4 equal segments we have either a square or rhombuses of arbitrarily low area, with the maximum to minimum ratio unbounded.

[3] gives a very pleasing 3D example which answers the above basic question in the negative. Consider the two objects - the rhombicuboctahedron and the pseudorhombicuboctahedron. Both these convex polyhedrons have the same sets of faces - 8 triangles, 18 squares. So if one is given only the faces, one is not sure as to the polyhedron from which the faces were broken off. Indeed, both polyhedrons are solutions to the reconstruction problem and they differ by a twist.

*A simpler example:* Consider a cube of side a and 2 square pyramids of base of side a and height less than a/2. Two different convex polyhedrons with same face sets result by attaching the pyramids on (1) opposite faces of the cube and (2) adjacent faces of the cube.

**Question:** What are the non-trivial invariant quantities among all distinct convex polyhedrons formable from any given face set?

**Proposition:** an obvious and very tempting candidate is volume. In both examples given above, the volume is same for convex polyhedrons formed from the same face set.

**Counterexample:** Sets of polygonal faces exist that can form different convex polyhedrons with different volumes. Indeed, [4] describes three convex objects – the tetracontrahedron, the icosagonal dipyramid and the decagonal dipyramidal antiprism. The latter two clearly have different volumes. All result from the same set of 40 identical isosceles triangular faces.

**Further questions:**

1. Given a set of n disconnected faces, let a certain number of distinct convex polyhedrons be assembled. Can some bound on their number, in terms of the number of faces be found? And if a face set were found that yields different convex polyhedrons of different volumes, is there a bound on the ratios between the max and min volumes? Intuitively, in the large n limit (infinitely many faces), this volume ratio could be arbitrarily large; but for finite n, we guess this ratio has a finite bound (unlike in 2D where, as noted above, the area ratio can be unbounded even for convex polygons assembled from the same finite edge set).

2. Are there face sets where not all faces are identical (or even with all faces different) and which generate distinct polyhedrons with different volumes? The example in [4] has all faces identical.

3. Can one have an upper bound on the total different shapes among the polygons in any set of faces that can yield convex polyhedrons of different volumes - and in face sets which yield polyhedrons with different shapes with maybe the same volume?

4. What happens in higher dimensions - convex polytopes assembled from a set of hyperfaces?

## 3. On 'Rectangling Rectangles'

In this section, we present 2 questions on rectangling rectangles – formation of rectangular layouts putting together smaller rectangles. The questions derive from the basic question presented at [5] and discussed in detail in [6]. The phrase 'rectangling the rectangle' was introduced by Michael Brand ([6]). We deal with only 2D versions of the questions – higher dimensional versions are straightforward to state.

*1. Find a set of rectangles all of same perimeter (let us call such rectangles isoperimetric rectangles) but with different areas which together form a neat big rectangle.*

*Note:* the problem of forming such a rectangular layout of rectangles with same area but different perimeters was attacked long ago (see [7]).

*A simple experimental answer:* The following 7 isoperimetric rectangles (their x and y dimensions listed) together form a big 24 X 18 rectangle:

10.000000, 9.500000
16.000000, 3.500000
6.000000, 13.500000
15.500000, 4.000000
5.500000, 14.000000
18.500000, 1.000000
2.500000, 17.000000.

The above set of 7 isoperimetric tiles, when arranged in an expanding spiral with the first one at the center, gives a big rectangular layout. Our experiments show that infinitely many such spiral layouts exist - all containing 7, 8 and 9 tiles with the dimensions of each rectangular tile rational.

**Conjectures:** It appears that no such rectangular layouts can be made with less than 7 isoperimetric tiles (we mean any rectangular layout, not only spirally arranged ones). There may be no layouts where (1) the big rectangle is a square (at least when the dimensions of the tiles are rational) or (2) with the number of isoperimetric tiles arbitrarily large (indeed, we know of no rectangular layout formed with more than 10 isoperimetric tiles). We do not yet know of any solution where the ratios between dimensions of the tiles are irrational.

*Note:* Similar questions could be asked with other shapes than rectangles. For example, any triangle appears to readily allow partition into any number of equal area triangles with different perimeters or into isoperimetric triangles with different areas. One can also consider a triangle partitioned into quadrilaterals and vice versa.

*2. Given a number n, find that set of n rectangular tiles of any area and perimeter (of course, one could choose some of these tiles as identical, some as squares, whatever) which together form the largest possible number of rectangular layouts with different perimeters. Obviously, every rectangular layout should use all rectangular tiles in the chosen set and the tiles should be arranged without holes or overlaps. Equivalently, one could try, given a number m, to find a set of n rectangular tiles which give exactly m layouts with n being the minimum possible.*

*Simple Example:* Let n =4. Choosing 4 unit squares as the tile set will give 2 different rectangular layouts - the 2X2 square and the 4X1 rectangle. But, if we choose a set of two 4x1 rectangles and two 2x1 rectangles, we can form three different rectangular layouts - a 12x1, a 4x3 and a 6x2. That is clearly better than the 4 unit squares.

We denote by n, the number of tiles to be used and with $L(n)$, the maximum number of rectangular layouts possible with n tiles (with the tiles chosen in the best possible way). It is easy to see that $L(n+1)$ is not less than $L(n)$ for any n. Indeed, consider the best set of n rectangles which give $L(n)$ layouts and divide any one rectangle in this set into 2 rectangles; now there are n+1 tiles and all arrangements that give $L(n)$ can from these. Experiments indicate $L(n)$ is a slowly growing function so obviously, it would not always increase whenever n is incremented by 1.

We now suggest a way to find rectangular tile sets with n members that can give large numbers of different rectangular layouts. We do not claim that this approach is optimal and gives $L(n)$ accurately.

A number that has more distinct factors than all numbers below it is said to be a highly composite number (hcn) [14].

The first few elements in the series of hcns:1, 2, 4, 6, 12, 24, 36, 48, 60, 120, 180, 240, 360, 720, 840, 1260, 1680, 2520, 5040, 7560, 10080, 15120, 20160, 25200, 27720, 45360, 50400, 55440, 83160, 110880, 166320, 221760, 277200, 332640, 498960, 554400, 665280, 720720, 1081080, 1441440, 2162160,…        (series 1)

The first few elements in the series of the number of factors of hcns:1, 2, 3, 4, 6, 8, 9, 10, 12, 16, 18, 20, 24, 30, 32, 36, 40, 48, 60, 64, 72, 80, 84, 90, 96, 100, 108, 120, 128, 144, 160, 168, 180, 192, 200, 216, 224, 240, 256, 288, 320, 336, 360, 384, 400, 432, 448, 480, 504, 512, 576, 600, 640, 672, 720, 768, 800, 864, 896,…        (series 2)

Both the above series are given in [15].

Consider any hcn, say 60. It has 12 factors. Consider a set of 60 identical rectangles of width 1 unit and some large length L. This set of tiles with combined width 60 will readily give rectangular layouts of widths equal to each of the 12 factors of 60 – from 1 to 60. Thus with 60 rectangles, we have 12 different layouts.

However, based on an hcn greater than 60, we can find another set of 60 rectangles which achieve more than 12 different rectangular layouts. Indeed, consider the set: { 20 1XL rectangles, 20 2XL rectangles, 20 3XL rectangles}. The combined width of all these rectangles is 120, an hcn larger than 60. With the 60 tiles in the new set, we can easily form rectangular layouts the widths of which could be any of the factors of 120 except 1 and 2 – a total of 14 different widths and layouts. We now explain this approach a little more:

The series of numbers which are the sums of n successive integers from 1: 1, 2, 6, 10, 15, 21, 28, 36, 45, 55, 66, 78, 91, 105, 120, 136, 153, 171, 190, 231, 253, 276, 300, …                                                                                  (series 3)

If m is the i'th element in series 3 (i.e. m is the sum of all positive integers up to i), and h be an hcn divisible by m and let d be the value h/m. We observe: with the set of id rectangles all with the same large length and widths chosen as {d rectangles of width 1, d rectangles of width 2,…, d rectangles of width i), we can form layouts of width equal to each of the factors of h greater than or equal to i.

*Example:* h = 60; m = 15; i = 5; d = 4. A rectangle set with 20 elements (4 rectangles each of width 1,2,3,4,5) can form layouts all of area 60 and with widths 5, 6, 10, 12, 15, 20, 30 and 60, i. e all factors of 60 greater than or equal to 5. Thus with 20 tiles, we have 8 layouts. Thus our suggestion is essentially to consider the series 1 and 3 and to choose numbers suitably.

This method cannot apply to all possible numbers of tiles/layouts because series 2 of the number of factors of hcns covers only a small fraction of all numbers. But we observe that for every prime, there is a (sufficiently large) hcn for which it is a factor; so even for larger elements in series 3, hcns can be found that are divisible by that element. Further, if the method is applicable to some hcn we could extend it to some more numbers as follows:

*Example:* Consider the above h =60, m=15, i=5, d=4 example that gives 8 layouts with 20 tiles. Let all tiles have the length L = 118 units. Cut only one of the 1X118 rectangles in tile set into 2 1X59 pieces. Now, there are 21 rectangular tiles and they easily give 9 layouts – all 8 layouts possible earlier plus a layout of width 59 and length 120.

**Remarks:** We have not examined layouts of rectangular tiles with irrational ratios between their dimensions. We could also consider sets of rectangles with varied dimensions (with as many of their lengths and widths different from one another as possible) and which give many rectangular layouts. Indeed, the previous problem of isoperimetric tiles forming a rectangle is the special case where every tile dimension is unique but with the constraint that their perimeters have to be the same – we have not found such sets of tiles which give *more than one* rectangular layout. We also do not know if requiring all dimensions plus all perimeters (and in addition, areas) of tiles to be different will yield tile sets which can give *multiple* rectangular layouts. Of course, if all tile dimensions plus areas and perimeters have to be unique, we can easily form layouts which are squares (in isoperimetric case, as noted above, we have no such example yet): indeed, a spiral layout such as the one described above for the isoperimetric case has all tile dimensions and tile areas unique; if this layout is scaled in x or y direction into a square, the perimeter of each tile will also become unique.

We could also consider the same basic question of this section for other shapes: for instance, sets of triangles which add up to multiple quadrilaterals (for which a possible – but not necessarily optimal - answer could be: sets of right triangles which pair up into rectangles which in turn form multiple rectangular layouts as discussed above), sets of triangles which form multiple triangles and so forth. We conclude this section noting that for the method proposed above, we have no proof of optimality even for the cases it appears to work well.

## 4. Maximizing and Minimizing the Diameter

**Question:** In 2D, if both area and perimeter are specified, which convex shape has maximum diameter? And which convex shape has the minimum diameter for given area and perimeter? The generalizations of these questions to higher dimensions are straightforward to state (Note: Diameter and Minimum width are defined in problem 1 above).

There are two obvious ways to approach this question. We could (1) fix area at some value A and vary only the required perimeter or (2) fix the perimeter at some value p and vary the area.

*Note:* This question is different from the extensively studied isodiametric and isoperimetric problems for convex polygons and regions - for example, see Messinghof ([8]).

With area A given, the perimeter has to be more than a critical value for a convex shape to be even possible. And at that critical perimeter, the only shape possible is a circle (obviously it gives both max and min diameter).

*Remarks:* In 2D, there are six possible questions: Specify 2 of the three quantities {diameter, area, perimeter} and then maximize or minimize the 3rd quantity; some of these questions have been answered earlier (for example, see [8], [11]) but to our knowledge, not all. Generalizations of these questions to higher dimensions appear obvious.

We now discuss some partial answers to the maximum-minimum question on diameter.

**Minimum diameter in 2D:** The mean width w of a convex figure is always proportional to its perimeter, p. Indeed, p = π w (see [9]). This implies: Keeping perimeter fixed, as the specified area is varied (within a range as shown below), the convex figure that minimizes the diameter is always a figure of constant width, provided a figure of constant width exists for that area and perimeter ([12]). With perimeter p, there is no convex region possible with area greater than that of the circle with that perimeter. Among all figures of constant width of 1 unit and perimeter π, the circle has greatest area (π/4 units) and the Reuleaux triangle the least area (1/2(π - √3) = 0.7 units approximately) – see [10]. With perimeter fixed at the some value p, if the required area is reduced below the corresponding Reuleaux area, the diameter of the required convex figure will be more than that of Reueaux. Indeed, because perimeter of the convex figure has to remain constant, its average width (=p/π) has to be constant. The figure will not be of constant width for area less than Reuleaux; and since the diameter is the maximum value of width, the latter has to increase.

**Lemma:** As noted above, among all constant width convex 2D regions of fixed perimeter, the circle has largest and Reuleaux the smallest area. It is also true that for every intermediate area (and with the same perimeter), there exist some convex region(s) of the same constant width.

**Proof:** if K and L are convex bodies of constant width 1, then (1-t)K + tL is a convex body of constant width 1 for every t in [0, 1]. Since the area of (1-t)K + tL varies continuously, it spans a segment.

We now consider: What shape will minimize diameter if required area is less than that of the Reuleaux triangle with perimeter kept constant at π units?

**Conjectures:** the minimal diameter convex body for the given fixed perimeter and area less than a critical value C of approximately 0.57 units (C is obviously less than 0.7, the approximate area of a Reuleaux triangle with diameter 1 unit) is probably, the sector of a circle. When the required area is exactly C, this sector has diameter nearly 1.045 (this diameter is actually the *radius* of the circle which contains the sector!) and angle measure π/3. For further reduction of required area with perimeter kept constant, the best shape remains a sector and its angle measure continuously reduces and diameter increases. Between the areas 0.7 (the Reuleaux area) and C, when the required area is continuously reduced, the minimum diameter shape 'morphs' (we do not know the precise sequence of shapes it goes thru) from Reuleaux to the sector and diameter increases slowly from 1 to 1.045.

**Maximum diameter in 2D:**

**Claim:** The 2D convex region that maximizes the diameter for specified area and perimeter is the intersection of two congruent disks with appropriate diameters and distance between centers (how a convex lens is diagrammed). In 3D, the convex shape which maximizes diameter with volume and surface area specified could be a spindle – the 3D analog of the lens.

**Explanation:** Sheng Liang et al [11] prove the following 'final theorem': the convex 2D region fully containing a given line segment and having a specified perimeter and the largest possible area is a symmetric lens formed by two identical arcs and the given segment as diameter. The above quoted answer follows readily. Indeed, if some other convex shape S gave the highest diameter, then by this theorem, a lens exists with the same diameter segment as S and with same area as S but with less perimeter. And obviously, for any such lens L, another lens with same area as L and greater perimeter and hence larger diameter than L exists.

*Note:* Considering slight deformations of the convex lens, we see that there exist in its neighborhood, convex shapes with same area and perimeter but less diameter and also lower minimum width. So, with area and perimeter fixed, maximizing the diameter does not automatically minimize the minimum width. It appears that for fixed area and perimeter, the least minimum width is given by a rectangle and the highest min width by an isosceles triangle.

**Acknowledgements:** This paper is partially based on several posts made at nandacumar.blogspot.com during the period September 2012 – February 2013. As was acknowledged there, the answers to the problems on maximizing and minimizing the diameter were mostly given by Prof. Roman Karasev. Thanks to N. Ramana Rao for contributing to the polyhedron reconstruction question.